\documentclass{article}

\usepackage{fullpage} 
\usepackage{amssymb} 
\usepackage{graphics}
\usepackage{epsfig}

\def\ds{\displaystyle}
\def\Rset{{\rm I}\!{\rm R}}
\def\qed{\hfill \ensuremath{\Box}}

\begin{document}

\title{Effects of spatial structure and diffusion\\
on the performances of the chemostat}
\author{{\sc I. Haidar and A. Rapaport}\\[2mm]
UMR INRA/SupAgro 'MISTEA' and EPI 'MERE'\\
2, pl. Viala 34060 Montpellier, France\\
e-mails: {\tt ihabhaidar@yahoo.com,rapaport@supagro.inra.fr}}
\date{}

\maketitle

\noindent {\em Abstract.} Given hydric capacity and nutrient flow of a chemostat-like system, we analyse the influence of a spatial structure on the output concentrations at steady-state. Three configurations are compared: perfectly-mixed, serial and parallel with diffusion rate. We show the existence of a threshold on the input concentration of nutrient for which the benefits of the serial and parallel configurations over the perfectly-mixed one are reversed. In addition, we show that the dependency of the output concentrations on the diffusion rate can be non-monotonic, and give precise conditions for the 
diffusion effect to be advantageous.
The study encompasses dead-zone models.

\section{Introduction}
The chemostat is a popular apparatus, invented simultaneously by Monod \cite{Monod} and Novick \& Szilard \cite{Novick}, for the so-called {\em continuous culture} of micro-organisms. It has the advantage to study bacteria growth at steady state, in contrast to batch cultivation. In the classical experiments, the medium is assumed to be perfectly mixed, that justifies mathematical models described by systems of ordinary differential equations \cite{Smith}. 
The chemostat model is also used in ecology for studying populations of micro-organisms, such as lake plankton or wetlands ecosystems.
In natural ecosytems, or in industrial applications that use large bioreactors, the assumption of perfectly mixed medium is questionable. This is why spatial considerations have been introduced in the classical model of the chemostat, such as the gradostat model \cite{Wimpenny} that is a series of interconnected chemostats (of identicalvolumes).
Segregated habitats are also considered in lakes, where the bottom 
can be modeled as a {\em dead zone} and nutrient mixing between the two zones is achieved by diffusion rate \cite{Nakaoka}. The consideration of dead zones is also often used in bioprocesses modelling \cite{Levenspiel,Hu,Grobicki,Roux,Roca,Valdes,Saddoud}.

Series of chemostats, instead of single chemostat, have shown to potentially improve the performances of bioprocesses, reducing the total residence time \cite{Hill,Luyben,Harmand03,Harmand04,Harmand05}
or allowing species persistence \cite{Stephanopoulos,Rapaport}. These properties have of course
economical impacts for the biotechnological industry, and 
there is a significant literature on the design of series of reactors and comparison with plug-flow reactors (that can be seen as the {\em limiting case} of an arbitrary large number of tanks of arbitrary small volumes) \cite{Riet,Tramper,Nelson,deGo,Doran,Drame1,Drame2}. 
Sometimes a radial diffusion is also considered in plug-flow reactors \cite{Grady}, but surprisingly, configurations of tanks in parallel have been much less investigated \cite{Levenspiel}. One can argue that knowing input rates and volumes of tanks in parallel, there dynamical characteristics can be studied separately, and there is no need of devoting a specific study for these configurations.
 This is no longer the case if one considers a passive communication 
between the tanks, through a membrane for instance. In saturated soils or wetlands, a spatial structure could be 
simply represented by separated domains with diffusive communication.
This consideration is similar to patches models or {\em islands models}, commonly used in ecology \cite{MacArthur,Hanski}, or lattice differential equations \cite{SV09}.
For instance, a recent investigation studies the influence of such structures 
on a consumer/resource model \cite{Gravel}. Consumer/resource models in ecology are similar to chemostat models, apart the source terms that are modeled as constant intakes of nutrient, instead of dilution rates that one rather met in liquid media. 

In this paper, we propose to bring
new insight on parallel configurations of chemostats with communication, in a spirit different than the one usually taken 
in bioprocesses design.  One usually chooses a target for the output concentration of substrate, and looks for minimizing the total volume, or equivalently the residence time, among all the configurations that provide the same desired output at steady state.
Here, we fix both the total hydric volume and the input flow and study the
input-output map at steady-state, investigating the role of 
the spatial structure on the performances of the system.
The performance is here measured by the level of 
substrate that is degraded by the system, and collected at the output. 
We draw precise comparisons
between the three configurations: perfectly mixed, serial and parallel (with diffusion rate) with the same total hydric volume and flow rate.
This set of configurations is far to be exhaustive, being
limited to two compartments only, but it is a first attempt 
to grasp this input-output map of a structured chemostat,
and study how a spatial structure can modify this map, and
what are the key parameters.
We believe that this study is of interest for the modelling of ecosystems such as saturated soils for which it is not easy to know the 
spatial structure, and where one has only access to {\em input-output} observations of the substrate degradation.\\

The paper is organized as follows. In Section \ref{section-model}, we present the three configurations under investigation and 
give the equations of the models. The main part of the paper is devoted to the analysis
of the steady states, given in Section \ref{section-analysis}. The proofs of the
global stability of the equilibriums are postponed to the Appendix, for lightening the presentation. Finally, discussion and numerical simulations are given in Section \ref{section-simulation}. 

\section{The models}
\label{section-model}
\begin{center}
\begin{figure}[h]
\begin{center}
\includegraphics[height=7cm]{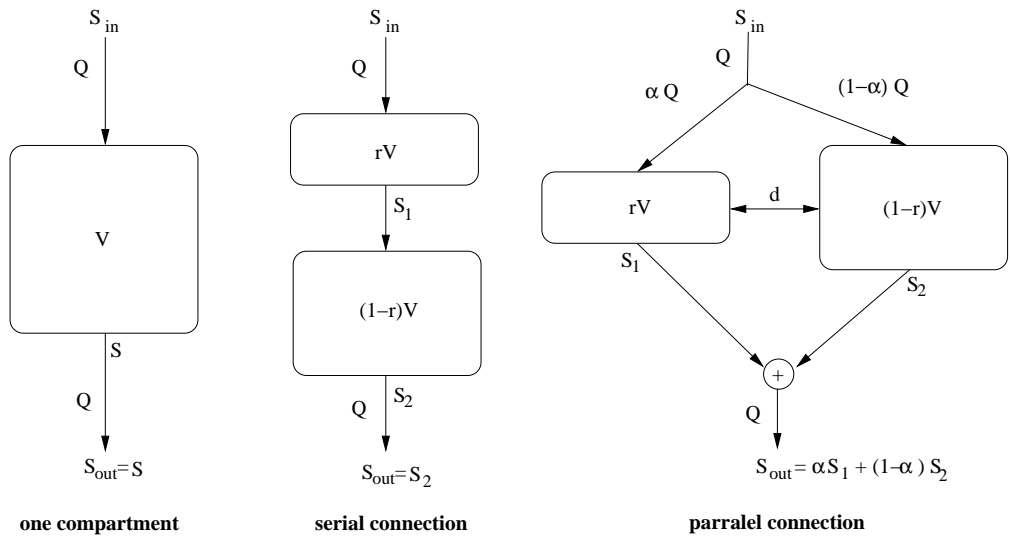}
\caption{The set of configurations under investigation.\label{fig-configs}}
\end{center}
\end{figure}
\end{center}

\begin{center}
\begin{figure}[h]
\begin{center}
\includegraphics[width=3cm]{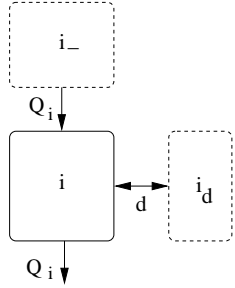}
\caption{Possible interconnections of a compartment.\label{fig-compartment}}
\end{center}
\end{figure}
\end{center}

The flow rate is labeled $Q$ and $V$ is the total capacity of the system. The three simple patterns we analyze are depicted on Figure \ref{fig-configs}: 
\begin{itemize}
\item[-] one single compartment of volume $V$
\item[-] two compartments of volume $V_{1}$, $V_{2}$ such that $V=V_{1}+V_{2}$ connected in cascade
\item[-] two compartments of volume $V_{1}$, $V_{2}$ such that $V=V_{1}+V_{2}$ in parallel with a diffusive connection.
\end{itemize}
We recall the dynamical equations of resource (nutrient) and biomass concentrations, respectively denoted by $S_{i}$ and $X_{i}$
in a compartment $i$ of volume $V_{i}$ fed from a compartment $i_{-}$ with a flow rate $Q_{i}$ and connected by diffusion rate $d$ to a compartment $i_{d}$ (see
Figure \ref{fig-compartment}).
\[
\begin{array}{lll}
\dot S_{i} & = & \ds -\frac{\mu(s_{i})}{y}X_{i} + \frac{Q_{i}}{V_{i}}(S_{i_{-}}-S_{i})
+ \frac{d}{V_{i}}(S_{i_{d}}-S_{i})\\
\dot X_{i} & = & \ds \mu(s_{i})X_{i} + \frac{Q_{i}}{V_{i}}(X_{i_{-}}-X_{i})
+ \frac{d}{V_{i}}(X_{i_{d}}-X_{i})
\end{array}
\]
For sake of simplicity of the analytical analysis, we assume that the growth function $\mu(\cdot)$ is a linear function of the resource concentration:
\[
\mu(S)=mS
\]
In Section \ref{section-simulation}, we shall consider Monod growth function 
and show that the qualitative results of our study are not changed.
The yield coefficient $y$ of the bio-conversion is kept equal to one (this is always possible by choosing the unit measuring the biomass).
It is convenient to write dimensionless concentrations:
for each concentration $C_{i}$  in the compartment $i$
($C_{i}$ can denote $S_{i}$ or $X_{i}$), we define
\[
c_{i} = m\frac{V}{Q}C_{i} \ ,
\]
and 
\[
r_{i}=\frac{V_{i}}{V} \ .
\]
We shall also consider that the time $t$ is measured in units such that
$Q=V$. Finally, we assume that the input concentration $S_{in}$
is large enough to avoid 
the (trivial) wash-out equilibrium to be the only steady-state 
in each compartment.

\section{Steady-state analysis of the three configurations}
\label{section-analysis}

\subsection{Configuration with one compartment}
The dynamical equations of the configuration with
a single compartment are
\[
\left\{\begin{array}{lll}
\dot s & = & -sx + s_{in}-s\\
\dot x & = & sx-x
\end{array}\right.
\]
The non-trivial equilibrium is $(1,s_{in}-1)$ under the condition $s_{in}>1$.
Then, one has
\[
s_{out}^{\star}=1 \ .
\]
{\em Remark.} This is a well known property from the theory of the chemostat that the output concentration at steady state is independent of the input concentration, provided this latter to be large enough
(i.e. $s_{in}\geq 1$).

\subsection{Serial connection of two compartments}
The dynamical equations of the model with two compartments in series (see Figure \ref{fig-configs}), assuming $r$ to be different to $0$ and $1$, are
\begin{equation}
\label{dyn2}
\left\{\begin{array}{lll}
\dot s_{1} & = & \ds -s_{1}x_{1} + \frac{1}{r}(s_{in}-s_{1})\\[2mm]
\dot x_{1} & = & \ds s_{1}x_{1}-\frac{1}{r}x_{1}\\[2mm]
\dot s_{2} & = & \ds -s_{2}x_{2} + \frac{1}{1-r}(s_{1}-s_{2})\\[2mm]
\dot x_{2} & = & \ds s_{2}x_{2}+\frac{1}{1-r}(x_{1}-x_{2})
\end{array}\right.
\end{equation}
with $r=V_{1}/V$.\\

\noindent {\bf Proposition 1.}
When $s_{in}>1/r$, there exists an unique equilibrium
$(s_{1}^{\star},x_{1}^{\star},s_{2}^{\star},x_{2}^{\star})$ of
(\ref{dyn2}) on the positive orthant. One has necessarily 
$s_{1}^{\star}=1/r$
and $s_{2}^{\star}<\min(1/r,1/(1-r))$. Furthermore, one has
\[
s_{out}^{\star} < 1 \Longleftrightarrow s_{in} > 1 + 1/r \ .
\]

\begin{center}
\begin{figure}[h]
\begin{center}
\includegraphics[height=4cm]{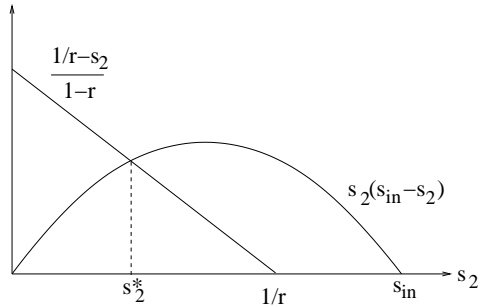}
\caption{Graphical determination of $s_{2}^{\star}$.\label{fig-dessin1}}
\end{center}
\end{figure}
\end{center}

{\em Proof.}
One can readily check that there exists a non-trivial equilibrium $(1/r,s_{in}-1/r)$ for the first compartment exactly when $s_{in}>1/r$.
Furthermore, this equilibrium is unique. Then, any equilibrium
for the overall system (\ref{dyn2}) has to be $(s_{2}^{\star},s_{in}-s_{2}^{\star})$ for the second compartment,
with $s_{2}^{\star}$ solution of the equation
\begin{equation}
\label{eq2}
s_{2}(s_{in}-s_{2})=\frac{1}{1-r}(1/r-s_{2})
\end{equation}
with $s_{2}^{\star}<1/r$. One can easily verify
that there exists a unique $s_{2}^{\star}$ solution of (\ref{eq2})
on $(0,1/r)$. Graphically, $s_{2}^{\star}$ is the abscissa 
of the intersection of the graphs (see Figure \ref{fig-dessin1})
of the polynomial function
\[
\phi(s_{2})=s_{2}(s_{in}-s_{2})
\]
and the affine function
\[
l(s_{2})= \frac{1}{1-r}(1/r-s_{2}) \ .
\]
Remark that $s_{in}>1/r$ implies the inequality $\phi(1/(1-r))>l(1/(1-r))$, from which one deduces $s_{2}^{\star}<1/(1-r)$.
Finally one can compare $s_{out}=s_{2}^{\star}$ 
with the value obtained in the configuration of one compartment:
\[
s_{out}^{\star} < 1 \Longleftrightarrow \phi(1)>l(1) \Longleftrightarrow s_{in} > 1 + 1/r \ . \qquad \qed
\]

\subsection{Parallel interconnection of two compartments}

The dynamical equations of the model with two compartments in parallel
and diffusion (see Figure \ref{fig-configs}), assuming $r$ to be different to $0$ and $1$, are the following
\begin{equation}
\label{dyn2d}
\left\{\begin{array}{lll}
\dot s_{1} & = & \ds -s_{1}x_{1} + \frac{\alpha}{r}(s_{in}-s_{1})
+\frac{d}{r}(s_{2}-s_{1})\\[2mm]
\dot x_{1} & = & \ds s_{1}x_{1}-\frac{\alpha}{r}x_{1}+\frac{d}{r}(x_{2}-x_{1})\\[2mm]
\dot s_{2} & = & \ds -s_{2}x_{2} + \frac{1-\alpha}{1-r}(s_{in}-s_{2})
+\frac{d}{1-r}(s_{1}-s_{2})\\[2mm]
\dot x_{2} & = & \ds s_{2}x_{2}-\frac{1-\alpha}{1-r}x_{2}+\frac{d}{1-r}(x_{1}-x_{2})
\end{array}\right.
\end{equation}
where the output concentration $s_{out}$ is given by
\[
s_{out} = \alpha s_{1} + (1-\alpha)s_{2} \ .
\]
The wash-out in both compartments corresponds to the trivial equilibrium
$(s_{in},0,s_{in},0)$, that leads to the trivial steady-state 
$s_{out}^{\star}=s_{in}$.
For convenience, we posit
\[
\alpha_{1}=\frac{\alpha}{r} \; , \;
\alpha_{2}=\frac{1-\alpha}{1-r} \ ,
\]
and assume, without any loss of generality that one has
$\alpha_{2} \geq \alpha_{1}$
(if it is not the case one can just exchange indexes $1$ and $2$).\\

\noindent {\em Remark.} One has necessarily $\alpha_{2}\geq 1$ and
$\alpha_{1}\leq 1$.\\

When $d=0$ (no diffusion), the equilibrium of the system 
can be determined independently in the two compartments as simple chemostats.
In this case, there is an unique globally stable equilibrium
$(s_{1}^{*},s_{in}-s_{1}^{\star},s_{2}^{*},s_{in}-s_{2}^{\star})$ in the non-negative orthant, where
$s_{i}^{*}=\min(\alpha_{i},s_{in})$ ($i=1, 2$).

When $d>0$, we define the functions
\[
\begin{array}{lll}
\phi_{2}(s_{1}) & = & \ds s_{1}+\frac{r}{d}(s_{in}-s_{1})(s_{1}-\alpha_{1}) \ ,\\[4mm]
\phi_{1}(s_{2}) & = & \ds s_{2}+\frac{1-r}{d}(s_{in}-s_{2})(s_{2}-\alpha_{2}) \ ,
\end{array}
\]
and 
\[
g(s_{1})=\phi_{1}(\phi_{2}(s_{1}))-s_{1} \ .
\]

\medskip

\noindent {\bf Proposition 2.} When $s_{in} > 1$ and $d>0$,
there exists a unique equilibrium $(s_{1}^{\star},x_{1}^{\star},
s_{2}^{\star},x_{2}^{\star})$ of (\ref{dyn2d}) in the positive orthant,
where ($s_{1}^{\star},s_{2}^{\star})$ is the unique solution of the system
\begin{equation}
\label{system}
s_{2}^{\star}=\phi_{2}(s_{1}^{\star}) \mbox{ and }
s_{1}^{\star}=\phi_{1}(s_{2}^{\star}) \ ,
\end{equation}
on the domain $(0,s_{in})\times (0,s_{in})$, with
$x_{i}^{\star}=s_{in}-s_{i}^{\star}$ ($i=1,2$).
Furthermore, $s_{1}^{\star}=s_{2}^{\star}=1$ when $\alpha_{2}=\alpha_{1}$
and 
\begin{equation}
\label{inequalities2d}
\alpha_{1} < s_{1}^{\star} <  
s_{2}^{\star} < \min(\alpha_{2},s_{in})
\end{equation}
when $\alpha_{2}>\alpha_{1}$.
\medskip

\begin{center}
\begin{figure}[h]
\begin{center}
\includegraphics[height=5cm]{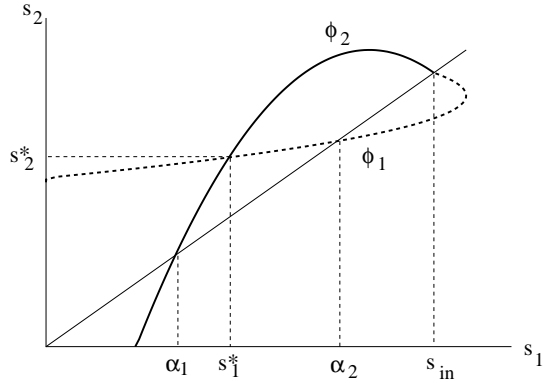}
\caption{Graphical determination of steady states
(when $\alpha_{1}<\alpha_{2}<s_{in}$).\label{fig-intersection}}
\end{center}
\end{figure}
\end{center}

{\em Proof.}
At equilibrium, one has
\[
\begin{array}{l}
r(\dot s_{1} +\dot x_{1})+(1-r)(\dot s_{2}+\dot x_{2})=0,\\
r(\dot s_{1} +\dot x_{1})=0,
\end{array}
\]
which amounts to write, from equations (\ref{dyn2d})
\[
\begin{array}{l}
\alpha(s_{in}-s_{1}^{\star}-x_{1}^{\star})+(1-\alpha)(s_{in}-s_{2}^{\star}-x_{2}^{\star})=0,\\
\alpha(s_{in}-s_{1}^{\star}-x_{1}^{\star})+d(s_{2}^{\star}+x_{2}^{\star}-s_{1}^{\star}-x_{1}^{\star})=0,
\end{array}
\]
or equivalently
\[
\underbrace{\left[\begin{array}{cc}
\alpha & 1-\alpha\\
\alpha+d & -d
\end{array}\right]}_{M}
\left(\begin{array}{c}
s_{in}-s_{1}^{\star}-x_{1}^{\star}\\
s_{in}-s_{12}^{\star}-x_{2}^{\star}
\end{array}\right)=
\left(\begin{array}{c} 0\\ 0 \end{array}\right)
\]
One has $det(M)=\alpha^{2}-\alpha-d\leq -d <0$ and deduces the property
\[
s_{1}^{\star}+x_{1}^{\star}=s_{2}^{\star}+x_{2}^{\star}=s_{in} \ .
\]
Consequently, an equilibrium in the positive orthant has to fulfill
$s_{i}^{\star} \in [0,s_{in}]$ for $i=1, 2$.
Replacing $x_{i}^{\star}$ by $s_{in}-s_{i}^{\star}$ in equations
(\ref{dyn2d}) at equilibrium, one obtains the equations
\begin{equation}
\label{equilibria2d}
\begin{array}{lll}
d(s_{2}^{\star}-s_{1}^{\star}) & = & 
    r(s_{in}-s_{1}^{\star})(s_{1}^{\star}-\alpha_{1})\\
d(s_{1}^{\star}-s_{2}^{\star}) & = & 
    (1-r)(s_{in}-s_{2}^{\star})(s_{2}^{\star}-\alpha_{2})
\end{array}
\end{equation}
which amounts to write that
$(s_{1}^{\star},s_{2}^{\star})$ is solution of the system
(\ref{system})
(see Figure \ref{fig-intersection})
or equivalently $s_{1}^{\star}$ is a zero of the function $g(\cdot)$.
When $\alpha_{2}=\alpha_{1}=1$, one can check that 
$s_{1}^{\star}=s_{2}^{\star}=1<s_{in}$ is solution of (\ref{system}).
When $\alpha_{2}>\alpha_{1}$, one has necessarily $\alpha_{1}<1$
and the condition $s_{in}>1$
implies $g(\alpha_{1})< 0$. We distinguish now two cases:
\begin{enumerate}

\item[] Case $\alpha_{2}< s_{in}$. If $\phi_{2}(\alpha_{2})\leq s_{in}$,
notice that one has $\phi_{2}(\alpha_{2})>\alpha_{2}$ and then
$g(\alpha_{2})>0$. If $\phi_{2}(\alpha_{2})> s_{in}$, notice that
$\phi_{2}(\alpha_{1})=\alpha_{1}<s_{in}$ and by the Mean Value Theorem, 
there exists $\tilde{s}_{2} \in (\alpha_{1},\alpha_{2})$ such that 
$\phi_{2}(\tilde{s}_{2})=s_{in}$ which implies 
$g(\tilde{s}_{2})=s_{in}-\tilde{s}_{2}>0$. In both cases, 
one deduces by the Mean Value Theorem the existence of
$s_{1}^{\star}\in(\alpha_{1},\alpha_{2})$ such that
$g(s_{1}^{\star})=0$.

\item[] Case $\alpha_{2}\geq s_{in}$. One has $g(s_{in})=0$ with 
\[
g^{\prime}(s_{in})= \frac{r(1-r)}{d^{2}}(\alpha_{1}-s_{in})(\alpha_{2}-s_{in})+\frac{1-s_{in}}{d}< 0 \ .
\]
Rolle and Mean Value Theorems allow to conclude the existence
of $s_{1}^{\star}\in(\alpha_{1},s_{in})$ such that
$g(s_{1}^{\star})=0$.
\end{enumerate}
In any case, we obtain the existence of $(s_{1}^{\star},s_{2}^{\star})$
solution of (\ref{system}) with
$s_{1}^{\star} \in (\alpha_{1},\min(\alpha_{2},s_{in}))$, 
that implies $s_{2}^{\star}=\phi_{2}(s_{1}^{\star})> s_{1}^{\star}$. 
But then 
$s_{1}^{\star}=\phi_{1}(s_{2}^{\star})< s_{2}^{\star}$ 
implies $s_{2}^{\star}<\min(\alpha_{2},s_{in})$.
Thus, the inequalities (\ref{inequalities2d}) are fulfilled.

Finally, notice that functions $\phi_{1}(\cdot)$, $\phi_{2}(\cdot)$ are 
both strictly concave, and steady states $(s_{1}^{\star},s_{2}^{\star})$ are
intersections of ${\cal G}_{1}$, 
the graph of the function $\phi_{1}(\cdot)$, and  ${\cal G}_{2}$ the
symmetric of the graph of $\phi_{2}(\cdot)$ with respect to the first 
diagonal.
Consequently, if $(s_{1}^{\star},s_{2}^{\star})$ is a steady state different 
from $(s_{in},s_{in})$,  ${\cal G}_{1}$ and ${\cal G}_{2}$ are
respectively above and below the line segment 
$(s_{1}^{\star},s_{2}^{\star})-(s_{in},s_{in})$. 
We conclude that there exists at most one non-trivial equilibrium.\qed\\

\noindent {\bf Corollary 2.} When $s_{in}>1$ and $d>0$, 
 the value $s_{1}^{\star}$ of the non trivial equilibrium
is the unique zero of the function $g(\cdot)$ on $(\alpha_{1},\min(\alpha_{2},s_{in}))$.
Furthermore, one has $g^{\prime}(s_{1}^{\star})> 0$.\\

\medskip

{\em Proof.}
When $\alpha_{1}=\alpha_{2}$, one has $s_{1}^{\star}=s_{2}^{\star}=1$ and one can easily check
\[
g^{\prime}(s_{1}^{\star})=\left(1+\frac{r}{d}(s_{in}-1)\right)\left(1+\frac{1-r}{d}(s_{in}-1)\right)-1 > 0 \ .
\]
When $\alpha_{2}>\alpha_{1}$, one has $g(\alpha_{1})<0$ and we recall
from the proof of former Proposition that $s_{1}^{\star}$ is the unique zero of $g(\cdot)$ on $(\alpha_{1},\min(\alpha_{2},s_{in}))$.
We conclude that $g$ is non decreasing at $s_{1}^{\star}$.
Notice that $\phi_{1}$ and $\phi_{2}$ are concave functions and
that
\[
\phi_{1}^{\prime}(\phi_{2}(s_{1}^{\star})))=1+\frac{1-r}{d}(s_{in}+\alpha_{2}-2s_{2}^{\star})>0
\]
implies
\[
g^{\prime\prime}(s_{1}^{\star})=\phi_{1}^{\prime\prime}(\phi_{2}(s_{1}^{\star})).\left[\phi_{2}^{\prime}(s_{1}^{\star})\right]^{2}+
\phi_{1}^{\prime}(\phi_{2}(s_{1}^{\star})).\phi_{2}^{\prime\prime}(s_{1}^{\star})<0
\]
We deduce that $g^{\prime}(s_{1}^{\star})$ cannot be equal to zero,
and consequently one has  $g^{\prime}(s_{1}^{\star})>0$. \qed \\

The global stability of the non-trivial equilibrium is proved in the Appendix.\\

\medskip
Proposition 2 defines properly the map
$d \mapsto s_{out}^{\star}=\alpha s_{1}^{\star}+(1-\alpha)s_{2}^{\star}$ for the unique non-trivial steady-state,
that we aim at studying as a function of $d$.
Accordingly to Proposition 2, $s_{out}^{\star}$ is equal to one for any 
value of the parameter $d$ in the non-generic case $\alpha_{2}=\alpha_{1}$.
We shall focus on the case $\alpha_{2}\neq\alpha_{1}$ (and without loss of generality we shall consider $\alpha_{2}>\alpha_{1}$).
We start by the two extreme situations: no diffusion and infinite diffusion.\\

\noindent {\bf Lemma 1.} For the non trivial equilibrium, one has
\[
s_{out}^{\star}(0) \geq 1 \Longleftrightarrow s_{in}\geq
s_{in}^{0}=\frac{r-\alpha^{2}}{r(1-\alpha)}
\]
with $s_{in}^{0} \in (1,2)$.\\

{\em Proof.}
Under the assumptions $s_{in}>1$ and $\alpha_{2}\geq \alpha_{1}$, we distinguish two cases when $d=0$.\\
If $s_{in}\geq\alpha_{2}$, one has $s_{1}^{\star}=\alpha_{1}$ and $s_{2}^{\star}=\alpha_{2}$. Then, one can write
\[
s_{out}^{\star}=\frac{\alpha^{2}}{r}+\frac{(1-\alpha)^{2}}{1-r}=1+\frac{(\alpha-r)^{2}}{r(1-r)}\geq 1 \ .
\]
If $s_{in} < \alpha_{2}$, on has $s_{1}^{\star}=\alpha_{1}$,
$s_{2}^{\star}=s_{in}$ and
\[
s_{out}^{\star}\geq 1 \Longleftrightarrow s_{in}\geq \frac{1-\alpha\alpha_{1}}{1-\alpha}=s_{in}^{0} \ .
\]
(recall that assuming $\alpha_{2}\geq \alpha_{1}$ imposes to have $\alpha<1$, and $s_{in}^{0}$ is well defined).
Notice that the number $s_{in}^{0}$ is necessarily larger than one because $\alpha_{1}\leq 1$, and one has also
\[
\alpha_{2}-s_{in}^{0}=\frac{(r-\alpha)^{2}}{r(1-r)(1-\alpha)}\geq 0.
\]
Consequently one concludes that $s_{out}^{\star}\geq 1$ exactly when 
$s_{in}\geq s_{in}^{0}$.
Finally, remark that one has
\[
s_{in}^{0}=\frac{r-\alpha^{2}}{r(1-\alpha)}
= 1-\frac{(\alpha-r)^{2}}{r(r-\alpha)}+\frac{r-\alpha}{1-\alpha}<2
 \ . \qquad \qed
\]

\noindent {\bf Lemma 2.} For $s_{in}>1$, the non trivial equilibrium fulfill
\[
\lim_{d\to+\infty} s_{1}^{\star}(d)=
\lim_{d\to+\infty} s_{2}^{\star}(d)=
\lim_{d\to+\infty} s_{out}^{\star}(d)=1 \ .
\]
{\em Proof.}
For any $d>0$, Proposition guarantees the existence of a unique non trivial equilibrium 
$(s_{1}^{\star},s_{2}^{\star}) \in (0,s_{in})\times(0,s_{in})$ 
that is solution of (\ref{equilibria2d}). 
When $d$ is arbitrary large, one obtains from (\ref{equilibria2d})
\[
\lim_{d\to+\infty} s_{1}^{\star}(d)-s_{2}^{\star}(d)=0 \ .
\]
From equations (\ref{equilibria2d}), one deduces also the 
following equality valid for any $d$
\[
r(s_{in}-s_{1}^{\star})(s_{1}^{\star}-\alpha_{1})+
(1-r)(s_{in}-s_{2}^{\star})(s_{2}^{\star}-\alpha_{2})=0 \ ,
\]
that can rewritten, taking into account the equality
$r\alpha_{1}+(1-r)\alpha_{2}=1$:
\[
(s_{in}-s_{1}^{\star})(s_{1}^{\star}-1)=(1-r)(s_{1}^{\star}-s_{2}^{\star})(s_{in}+\alpha_{2}-s_{1}^{\star}-s_{2}^{\star}) \ .
\]
Consequently, one has
\[
\lim_{d\to+\infty} s_{1}^{\star}(d)=
\lim_{d\to+\infty}s_{2}^{\star}(d)=1
\mbox{ or }
\lim_{d\to+\infty} s_{1}^{\star}(d)=
\lim_{d\to+\infty}s_{2}^{\star}(d)=s_{in} \ .
\]
If $\alpha_{2}<s_{in}$, the property $s_{1}^{\star}<\alpha_{2}$ 
valid for any $d>0$ implies that $s_{1}^{\star}$ cannot converges to
$s_{in}$.\\
If $\alpha_{2}\geq s_{in}$ and 
$\lim s_{1}^{\star}=\lim s_{2}^{\star}=s_{in}$, there exists $d$ such that 
$rs_{1}^{\star}+(1-r)s_{2}^{\star}>(s_{in}+1)/2$. Then, one has
\[
g^{\prime}(s_{1}^{\star})=
\frac{r(1-r)}{d^{2}}(s_{in}+\alpha_{1}-2s_{1}^{\star})(s_{in}+\alpha_{2}-2s_{2}^{\star})+\frac{s_{in}+1-2(rs_{1}^{\star}+(1-r)s_{2^{\star}})}{d}<0
\]
that contradicts Corollary 2. Finally, one has $\lim s_{1}^{\star}=\lim s_{2}^{\star}=1$ and consequently  $\lim s_{out}^{\star}=1$. \qed\\

We present now our main result concerning properties of the map
$d \mapsto s_{out}^{\star}(d)$ defined at the non-trivial steady-state.\\

\noindent {\bf Proposition 3.}
Assume $\alpha_{2}>\alpha_{1}$.
\begin{itemize}
\item[-] When $s_{in}\geq 2$, the map $d \mapsto s_{out}^{\star}(d)$
(for the non trivial equilibrium) is decreasing and $s_{out}^{\star}(d)>1$ for any $d\geq 0$.
\item[-] When $s_{in}< 2$,  the map $d \mapsto s_{out}^{\star}(d)$ 
(for the non trivial equilibrium)
admits a minimum in $d^{\star}<+\infty$, that is strictly less than one.
Furthermore, one has
\[
s_{in} > \underline s_{in}=\frac{2\alpha_{1}\alpha_{2}}{\alpha_{1}+\alpha_{2}}
\quad \Longrightarrow \quad d^{\star }> 0 
\]
with $\underline s_{in}<\min(2,\alpha_{2})$.
\end{itemize} 

\medskip

{\em Proof.}
Let differentiate with respect to $d$ the equations (\ref{equilibria2d}) 
at steady state:
\[
\begin{array}{lll}
\ds (s_{2}^{\star}-s_{1}^{\star})+d\left(\partial_{d}s_{2}^{\star}-\partial_{d}s_{1}^{\star}\right) & = & \ds \underbrace{r(s_{in}-2s_{1}^{\star}+\alpha_{1})}_{A}\partial_{d}s_{1}^{\star}\\
\ds (s_{1}^{\star}-s_{2}^{\star})+d\left(\partial_{d}s_{1}^{\star}-\partial_{d}s_{2}^{\star}\right) & = & \ds \underbrace{(1-r)(s_{in}-2s_{2}^{\star}+\alpha_{2})}_{B}\partial_{d}s_{2}^{\star}
\end{array}
\]
that can rewritten as follows
\[
\underbrace{\left[\begin{array}{cc} A+d & -d\\ d & -B-d
\end{array}\right]}_{\Gamma}
\left(\begin{array}{c} \partial_{d}s_{1}^{\star}\\
\partial_{d}s_{2}^{\star} \end{array}\right) = (s_{2}^{\star}-s_{1}^{\star})
\left(\begin{array}{c} 1 \\ 1 \end{array}\right)
\]
Remark that one has
\[
\begin{array}{lll}
A+d & = & d\phi_{2}^{\prime}(s_{1}^{\star})\\
B+d & = & d\phi_{1}^{\prime}(s_{2}^{\star})\\
det(\Gamma) & = & d^{2}(1-\phi_{1}^{\prime}(s_{2}^{\star})\phi_{2}^{\prime}(s_{1}^{\star}))=-d^{2}g^{\prime}(s_{1}^{\star})
\end{array}
\]
From the Corollary 2, one has $det(\Gamma)<0$ and one deduces that
the derivatives $\partial_{d}s_{1}^{\star}$,  $\partial_{d}s_{2}^{\star}$
are defined as follows
\begin{equation}
\label{partiald1s2}
\begin{array}{lll}
\partial_{d}s_{1}^{\star} & = & \ds (s_{2}^{\star}-s_{1}^{\star})\frac{-B}{det(\Gamma)}\\[4mm]
\partial_{d}s_{2}^{\star} & = & \ds (s_{2}^{\star}-s_{1}^{\star})\frac{A}{det(\Gamma)}
\end{array}
\end{equation}
Notice from inequalities (\ref{inequalities2d}) that we obtain $B>0$  and deduce 
$\partial_{d}s_{1}^{\star}>0$ for any $d$. With Lemma 2 we conclude that 
$s_{1}^{\star}(d)<1$ for any $d$.\\
From equations (\ref{partiald1s2}), we can write
\[
\partial_{d}s_{out}^{\star} = (s_{2}^{\star}-s_{1}^{\star})
\frac{\alpha B-(1-\alpha)A}{-det(\Gamma)}=
\underbrace{[\alpha_{1}(s_{in}-2s_{2}^{\star})-\alpha_{2}(s_{in}-2s_{1}^{\star})]}_{\sigma}
\frac{(s_{2}^{\star}-s_{1}^{\star})r(1-r)}{-det(\Gamma)}
\]
When $s_{in}\geq 2$, one has $A>0$ and then $\partial_{d}s_{2}^{\star}<0$.
With Lemma 2 we conclude that 
$s_{2}^{\star}(d)>1$ for any $d$. Then, one obtain the inequality 
\[
\sigma<
(s_{in}-2)(\alpha_{1}-\alpha_{2})\leq 0
\]
which proves with Lemma 2 that $s_{out}^{\star}$ is a decreasing function of d that converges to one.\\

When $s_{in}<2$, we write
\[
\sigma=(s_{in}-2)(\alpha_{1}-\alpha_{2})+2(\alpha_{1}(1-s_{2}^{\star})-
\alpha_{2}(1-s_{1}^{\star}))
\]
As $s_{1}^{\star}$ and $s_{2}^{\star}$ tend to one when $d$ takes arbitrary large values, we conclude that there exists $\bar d<+\infty$ such that
$\sigma>0$ for any $d>\bar d$ and consequently $s_{out}^{\star}$ is smaller than one and increasing for $d>\bar d$. We conclude 
that the map $d \mapsto s_{out}^{\star}(d)$ admits a minimum, say at $d^{\star}<+\infty$, that is strictly less than one.\\

When $d=0$, one has $s_{1}^{\star}=\alpha_{1}$ and $s_{2}^{\star}=\alpha_{2}$ if $s_{in}\geq \alpha_{2}$. Then, one obtains
$\sigma=s_{in}(\alpha_{1}-\alpha_{2})<0$.  
So the map $d \mapsto \partial_{d}s_{out}^{\star}(d)$ is decreasing
at $d=0$ and consequently $d^{\star}>0$.

When  $d=0$ with
$s_{in}< \alpha_{2}$, one has  $s_{2}^{\star}=s_{in}$  and then
$\sigma=2\alpha_{1}\alpha_{2}-s_{in}(\alpha_{1}+\alpha_{2})$,
for which we conclude
\[
\sigma<0 \Longleftrightarrow s_{in} > \frac{2\alpha_{1}\alpha_{2}}{\alpha_{1}+\alpha_{2}}=\underline s_{in} \ .
\]
Remark that this case is feasible because of the inequality
$2\alpha_{1}\alpha_{2}<\min(2,\alpha_{2})(\alpha_{1}+\alpha_{2})$.
We conclude that for $s_{in}$ larger than this last value, $d^{\star}$ is necessarily strictly positive. \qed\\

\noindent {\em Remark.} The particular case $\alpha=0$  corresponds 
to a configuration of a perfectly mixed tank of volume $(1-r)V$ connected to a {\em dead-zone} of volume $rV$. 
This is a way to approximate a non well-mixed tank or {\em segregated} bioreactors of total volume $V$, estimating the fraction of the volume occupied by the highly agitated area.

\section{Numerical computation and discussion}
\label{section-simulation}

Propositions 1 and 3 reveal the existence of a threshold 
on the value of the input concentration $s_{in}$
(equal to $2$ for our choice of the parameters units)
that reverses the performances of the serial and parallel configurations in terms of $s_{out}^{\star}$, compared to the single tank case (for which $s_{out}^{\star}=1$):
\begin{itemize}
\item[-] for $s_{in}>2$, there exist serial configurations 
such that $s^{\star}_{out}<1$ for $r$ large enough (i.e. the first tank has to be large enough), but any parallel configuration produces
$s^{\star}_{out}>1$,
\item[-] for $s_{in}<2$, there exists parallel configurations 
such that $s^{\star}_{out}<1$, while any serial configuration
has $s^{\star}_{out}>1$. There exists 
another threshold $s_{in}^{0} \in (1,2)$ such that
configurations with $s^{\star}_{out}<1$ require to have 
$d$ large enough when $s_{in}>s_{in}^{0}$ (cf Lemma 1).
\end{itemize}
Furthermore, the best performance of the parallel configuration is obtained
\begin{itemize}
\item[-] for arbitrary large values of $d$ when $s_{in}>2$ ,
\item[-] for a finite positive $d^{\star}$ when 
$s_{in}\in(\underline s_{in},2)$ (where the expression of
$\underline s_{in}$ is given in Proposition 3).
\end{itemize}

\begin{center}
\begin{figure}[h]
\begin{center}
\includegraphics[height=5.5cm]{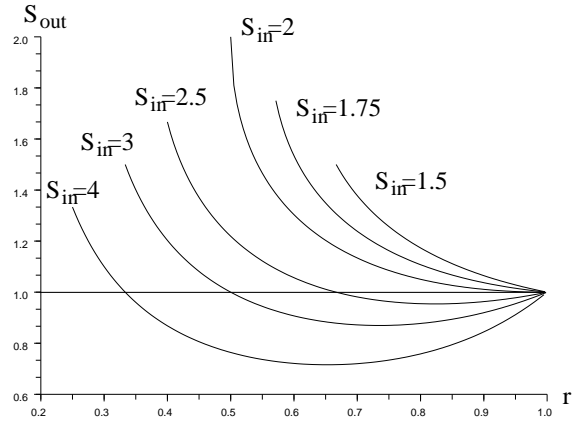}
\caption{Comparison of $s_{out}^{\star}$ for the serial configuration.\label{fig-comp2}}
\end{center}
\end{figure}
\end{center}

For the serial configuration, the graph of the function
$s^{\star}_{out}$ is plotted as function of $r\in [1/s_{in},1]$ 
on Figure \ref{fig-comp2} for different values of the input concentration $s_{in}$.

\begin{center}
\begin{figure}[h]
\begin{center}
\includegraphics[height=5.5cm]{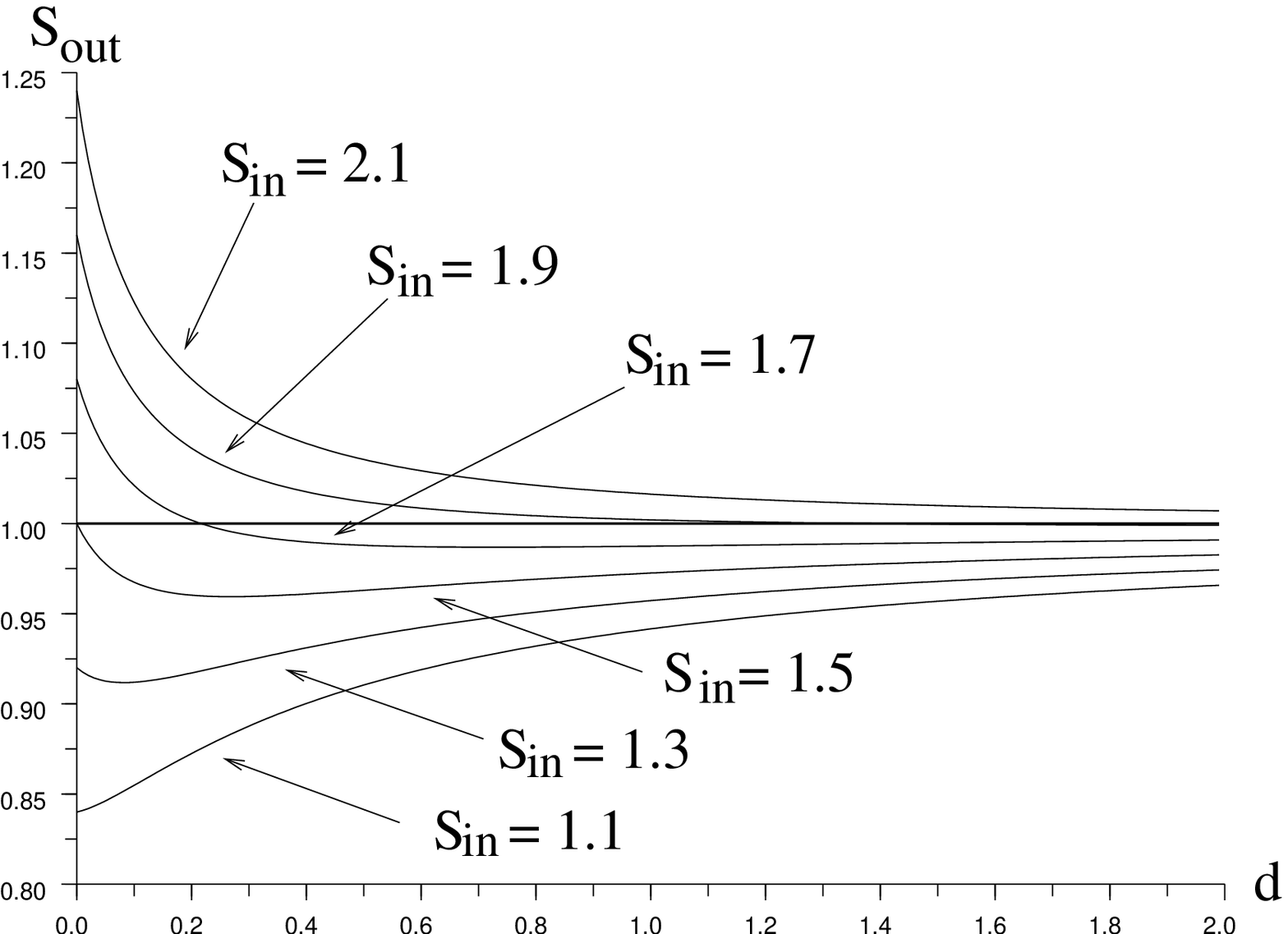}
\includegraphics[height=5.5cm]{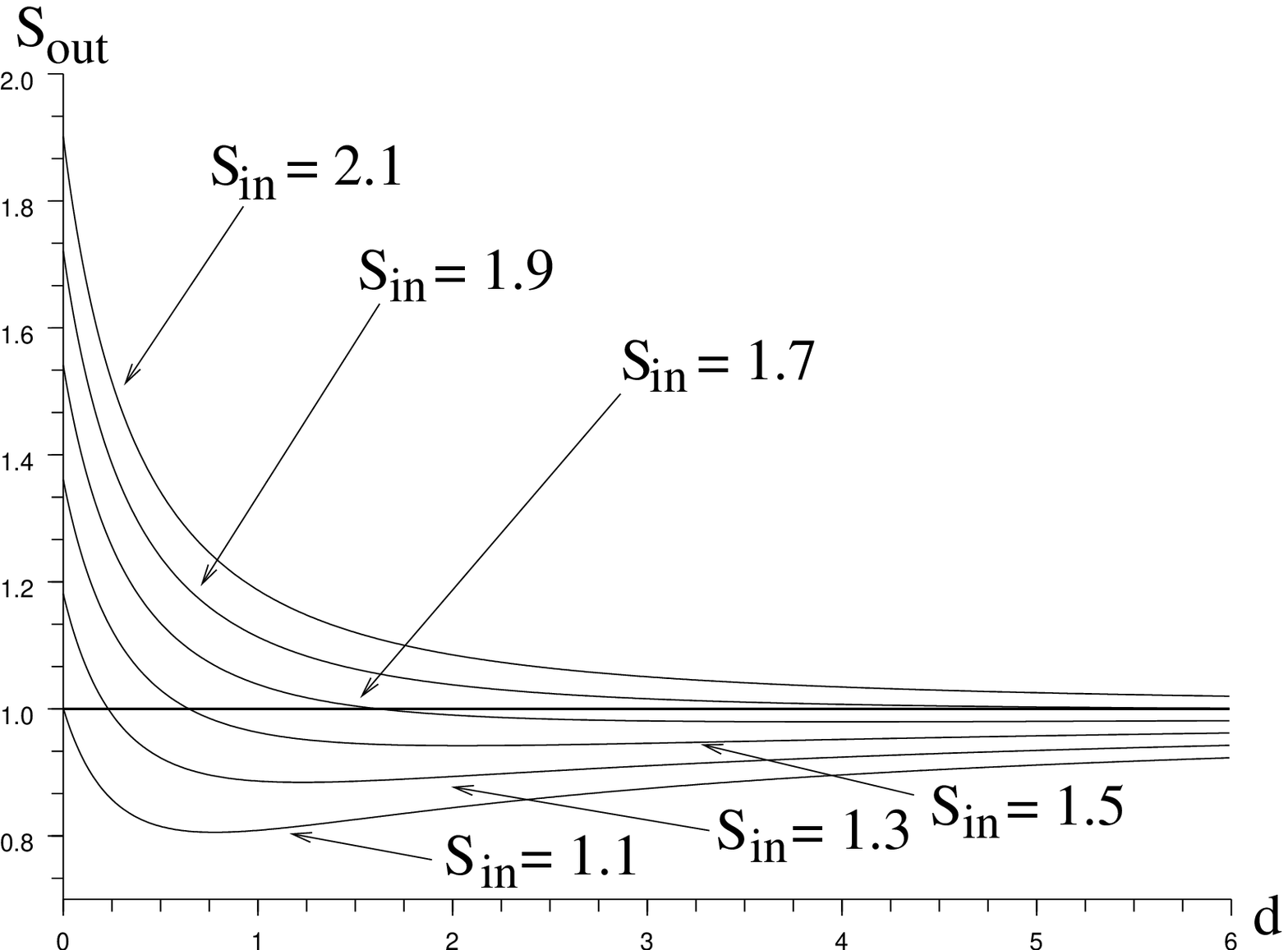}
\caption{Comparison of $s_{out}^{\star}$ for the parallel configuration ($\underline s_{in}>1$ on the left and $\underline s_{in}<1$ on the right)
.\label{fig-comp2d}}
\end{center}
\end{figure}
\end{center}

For the parallel interconnection, we depict on Figure
\ref{fig-comp2d} the two kind 
of configurations that occur, depending on whether the number
$\underline s_{in}$ is larger than one or not.
The values of the parameters are given on the table below
\begin{center}
\begin{tabular}{c|c|c||c|c|}
&$\alpha$ & $r$ & $\underline s_{in}$ & $s_{in}^{0}$\\
\hline
left figure & $0.6$ & $0.9$ & $1.14$ & $1.5$\\
right figure & $0.1$ & $0.9$ & $0.21$ & $1.09$\\
\end{tabular}
\end{center}

\bigskip

\begin{center}
\begin{figure}[h]
\begin{center}
\includegraphics[height=4cm]{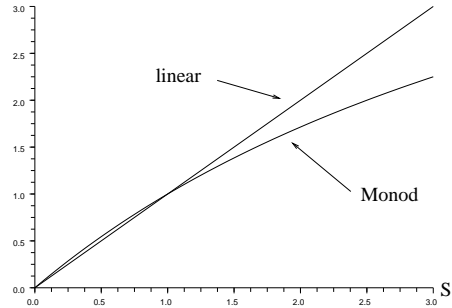}
\caption{Linear and Monod growth functions.
\label{fig-Monod}}
\end{center}
\end{figure}
\end{center}

The analytic analysis of Section \ref{section-analysis} has been conducted under the assumption of the linearity of the function $\mu(\cdot)$. 
It is often in microbiology that the growth rate $\mu(\cdot)$ presents a concavity, as described by the usual Monod (or Michaelis-Menten) 
function. 
We have computed numerically the same curves 
$s_{out}^{\star}(\cdot)$ than Figures \ref{fig-comp2} and 
\ref{fig-comp2d}, considering the Monod function
\[
\mu(S)=\frac{6S}{5+S}
\]
instead of the linear function (see Figure \ref{fig-Monod}).
This function has been chosen to fulfill  
$s_{out}^{\star}=1$ for the single tank configuration, 
guaranteeing the same steady state than the linear growth for
this configuration.

\begin{center}
\begin{figure}[h]
\begin{center}
\includegraphics[height=5.5cm]{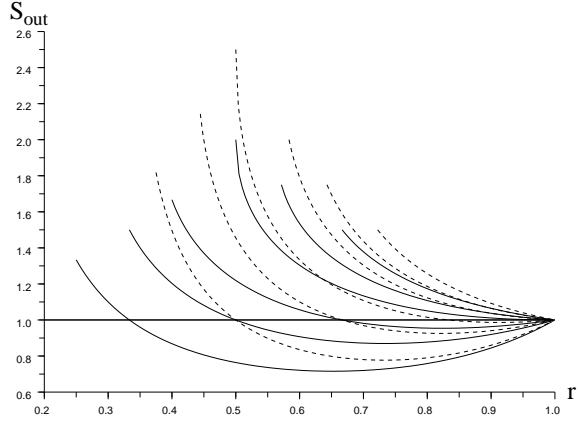}
\caption{Comparison of $s_{out}^{\star}$ for Monod (dashed) and linear (plain) for the serial configuration.
\label{fig-compare-serial}}
\end{center}
\end{figure}
\end{center}

\begin{center}
\begin{figure}[h]
\begin{center}
\includegraphics[height=5.5cm]{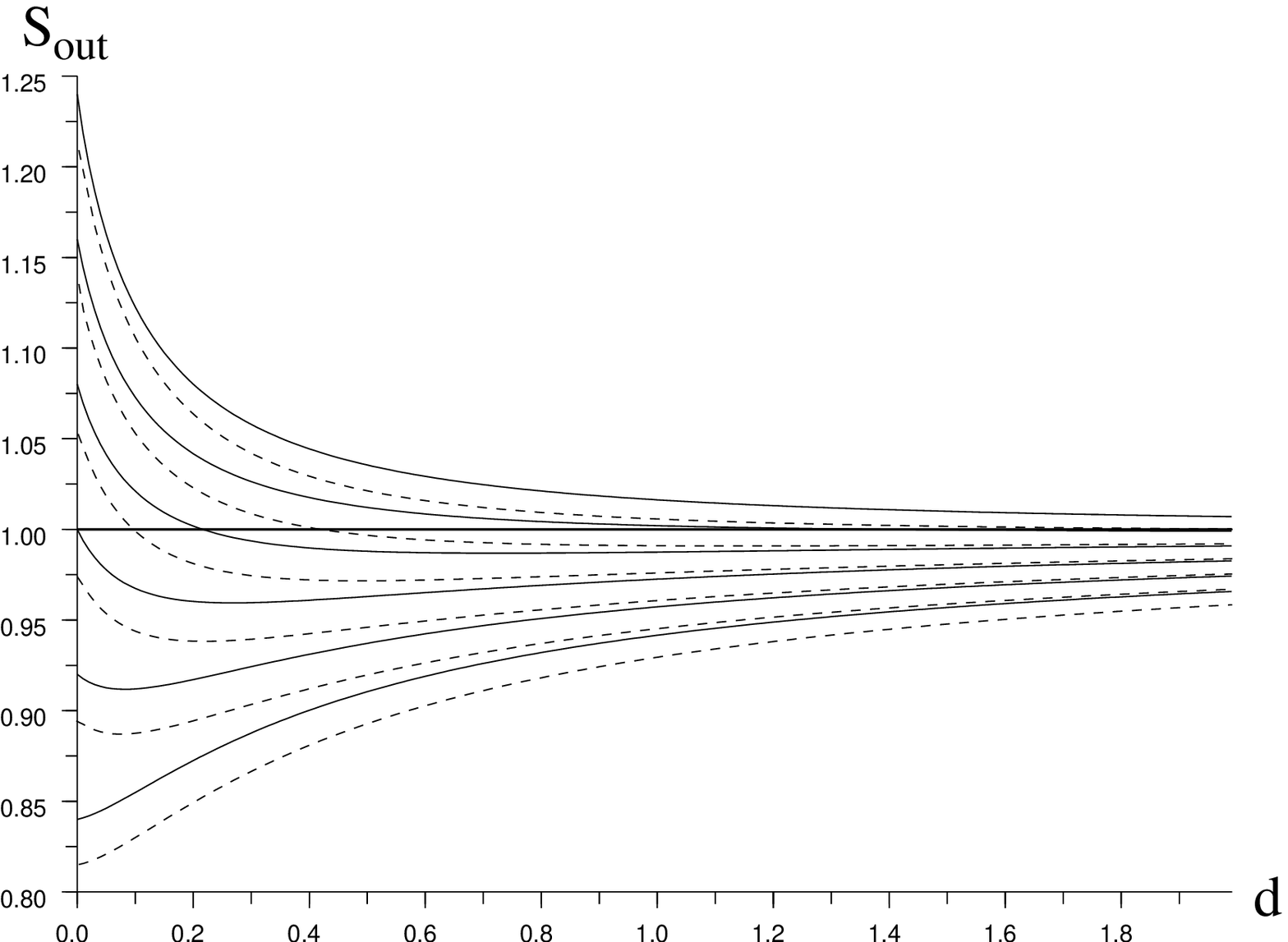}
\includegraphics[height=5.5cm]{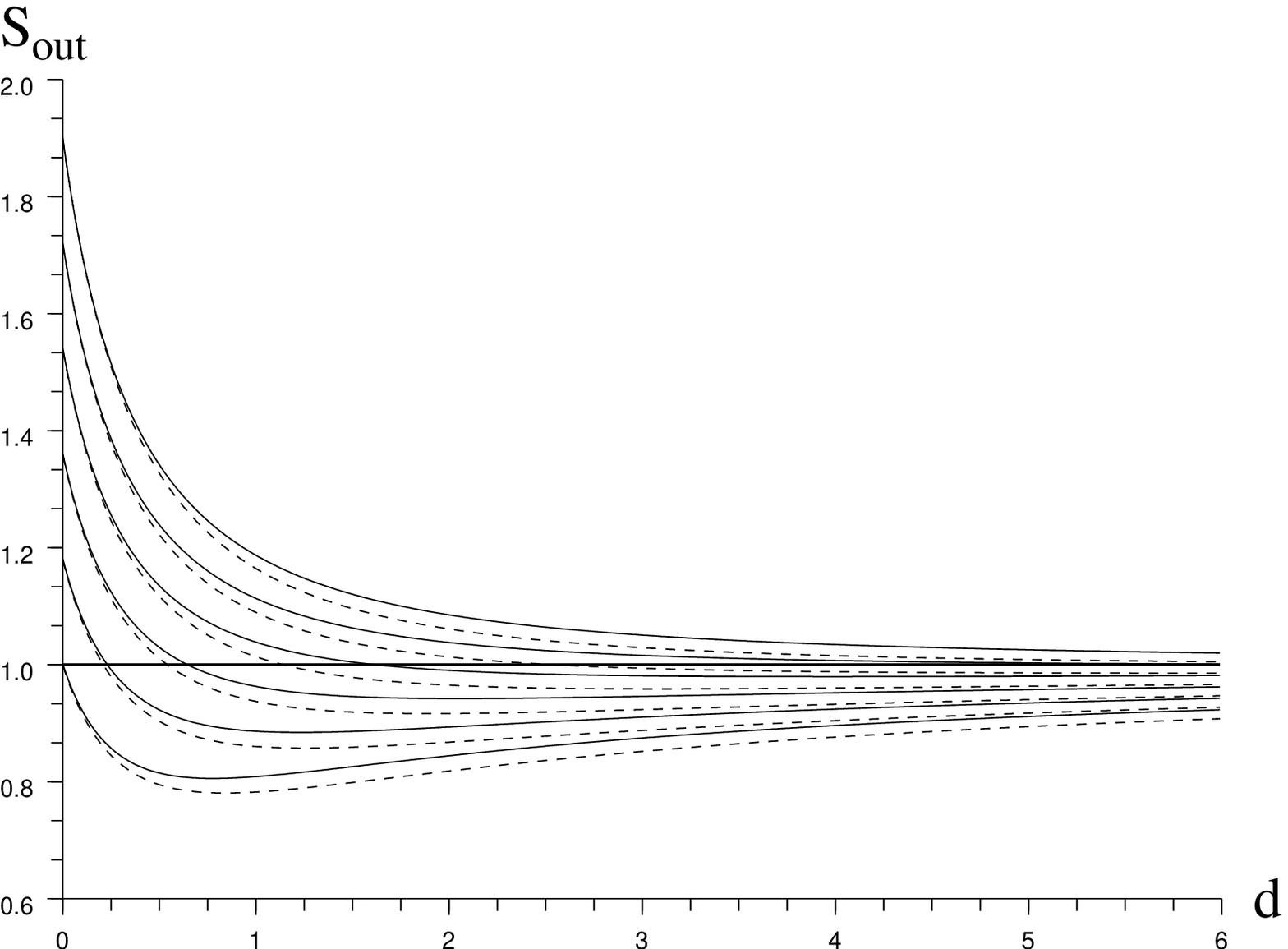}
\caption{Comparison of $s_{out}^{\star}$ for the parallel configuration ($\underline s_{in}>1$ on the left and $\underline s_{in}<1$ on the right)
.\label{fig-compare-parallel}}
\end{center}
\end{figure}
\end{center}

On Figures \ref{fig-compare-serial} and \ref{fig-compare-parallel}, we observe that the concavity of the growth function does not change qualitatively the theoretical results and the existence of threshold
for $s_{in}$ that favourites one of the configuration.
We notice on all the figures that the yield is better for the Monod function in the parallel configuration and worst for the serial one. This implies that the threshold on $s_{in}$, that was determined to
be equal to $2$ for the linear case is higher when the growth function is concave.\\

\noindent {\em Remarks.}
The serial configuration for the limiting value $r=1$ is equivalent to a single tank. This explains why all the curves on Figures
\ref{fig-comp2} and \ref{fig-compare-serial} coincide for this value of the parameter $r$.

For the parallel configuration with $\alpha=0.1$ and $r=0.9$ one has $\alpha_{2}=9$. This implies that for the limiting value $d=0$ the only equilibrium in the second tank is the wash-out when $s_{in}<9$. This is not the case for the first tank but the flow rate $\alpha Q$ being small, the output $s_{out}^{\star}$ remains closed to $s_{in}$ in any case, as one can see on Figures \ref{fig-comp2d} and \ref{fig-compare-parallel} for small values of the parameter $d$.

\section{Conclusion}

Given a flow rate and the total volume of a chemostat system, 
this study shows the existence of a threshold on the value of the input concentration $s_{in}$ such that above 
and below this threshold, serial and the parallel configurations are respectively the best ones with respect to the criterion of minimizing the output concentration $s_{out}^{\star}$ at steady state. For the parallel scheme, the best performances are obtained for a precise value of the diffusion parameter that is proved to be positive when $s_{in}$ is not too small. This study concerns also {\em dead-zone} configurations, as particular cases of the parallel configurations.

Whatever are the data of the problem, there always exists a configuration that is better than a single perfectly mixed tank. We have shown that the non-trivial steady states are unique and globally exponentially stable under the assumption of a linear increasing growth rate.

Finally, this study reveals the role of the structure of the space on the performances of simple ecosystems or bioprocesses.
The possibly non-monotonic influence of the diffusion parameter on the output steady state is not intuitive, and leave further investigations open for understanding or taking benefit of this property for natural ecosystems (such as saturated soils or wetlands) as well as for bioprocesses (such as waste-water treatments).
This result can be also of interest for reverse engineering when deciding which among serial or parallel configurations is better fit 
for the modeling of chemostat-like ecosystems, providing that one has an estimation of the hydric capacity of the system.\\

\noindent {\bf Acknowledgments.} This work has been achieved within the VITELBIO (VIRtual TELluric BIOreactors) program, sponsored by INRA and INRIA. The authors are grateful to this support.
The work is also part of the PhD thesis of the first author.

\section{Appendix: global exponential stability of the non-trivial equilibrium}

First, one can easily check that the domain $D=\Rset_{+}^{4}$ is 
invariant by the dynamics (\ref{dyn2}) and (\ref{dyn2d}).
We consider the 2-dimensional vector $z$ of variables $z_{i}=s_{in}-x_{i}-s_{i}$ ($i=1,2$) whose dynamics are respectively for the serial and parallel configurations
\[
\dot z =  A_{s}z=\left[\begin{array}{cc} -\frac{1}{r} & 0\\[2mm] 
\frac{1}{1-r} & -\frac{1}{1-r}\end{array}\right]z
\]
\[
\dot z = A_{p}z=\left[\begin{array}{cc} -\alpha_{1}-\frac{d}{r} & \frac{d}{r}\\[2mm] 
\frac{d}{1-r} & -\alpha_{2}-\frac{d}{1-r}\end{array}\right]z
\]
Notice that matrices $A_{s}$ and $A_{p}$ are Hurwitz~:
\[
\begin{array}{l}
\ds \mbox{tr}(A_{s}) = -\frac{1}{r}-\frac{1}{1-r}<0 \ , \quad
\mbox{det}(A_{s}) = \frac{1}{r(1-r)} >0\\
\ds \mbox{tr}(A_{p}) = -\alpha_{1}-\alpha_{2}-\frac{d}{r}-\frac{d}{1-r}<0 \ , \quad
\mbox{det}(A_{p}) = \alpha_{1}\alpha_{2} +\frac{d}{r(1-r)}>0
\end{array}
\]
So $z$ converges exponentially toward $0$ for both systems,
which implies that dynamics (\ref{dyn2}) and (\ref{dyn2d}) are dissipative, in the sense that any solution of (\ref{dyn2}) or (\ref{dyn2d})
in $D$ converge exponentially to the compact set
$K=\{(s_{1},x_{1},s_{2},x_{2})\in D \mbox{ s.t. }
x_{1}+s_{1}=s_{in} \mbox{ and } x_{2}+s_{2}=s_{in}\}$.\\

We recall a result from \cite[Theorem 1.8]{Thieme} 
that shall be useful in the following.\\

\noindent {\bf Theorem 4.}
Let $\Phi$ be an asymptotically autonomous semi-flow with limit semi-flow $\Theta$, and let the orbit ${\cal O}_{\Phi}(\tau,\xi)$ have compact closure. Then the $\omega$-limit set $\omega_{\Phi}(\tau,\xi)$ is non-empty, compact, connected, invariant and chain-recurrent by the semi-flow $\Theta$ and attracts $\Phi(t,\tau,\xi)$ when $t \to \infty$.

\subsection{The serial configuration}

\noindent {\bf Proposition  5.}
Under the condition $s_{in}>1/r$, any trajectory of
(\ref{dyn2}) with initial condition in $D$ such that $(s_{1}(0),x_{1}(0))\neq (s_{in},0)$ converges exponentially to the unique non-trivial steady-state $(s_{1}^{\star},x_{1}^{\star},s_{2}^{\star},x_{2}^{\star})$ given by Proposition 1.\\

{ \em Proof.}
Dynamics (\ref{dyn2}) has a cascade structure. It is straightforward to check that the solutions of the $(s_{1},x_{1})$ sub-system converges asymptotically towards the non-trivial equilibrium $(1/r,s_{in}-1/r)$ from any initial condition away from the wash-out equilibrium $(s_{in},0)$. 
From the convergence of $z_{2}$ toward $0$, we deduce that the $s_{2}$ variable has to converge to the bounded interval $[0,s_{in}]$ 
and that its dynamics can be written as a scalar non autonomous differential equation:
\begin{equation}
\label{nonautonomous}
\dot s_{2} = -s_{2}(s_{in}-s_{2}-z_{2}(t))+\frac{1}{1-r}(s_{1}(t)-s_{2})
\end{equation}
This last dynamics has the property to be asymptotically autonomous 
with the limiting differential equation:
\begin{equation}
\label{autonomous}
\dot s_{2} = f(s_{2})=-s_{2}(s_{in}-s_{2})+\frac{1}{1-r}(1/r-s_{2})
\end{equation}
Statement of Proposition 1 implies that this last scalar dynamics has a unique equilibrium $s_{2}^{\star}$ that belongs to $[0,s_{in}]$.
Furthermore, one has $f(0)>0$ and $f(s_{in})<0$. Consequently any solution of (\ref{autonomous}) in $[0,s_{in}]$ converges asymptotically to $s_{2}^{\star}$. Then applying Theorem 4, we conclude that any bounded solution of (\ref{nonautonomous}) converges to
$s_{2}^{\star}$. Finally any solutions of the $(s_{2},x_{2})$ sub-system converges asymptotically to $(s_{2}^{\star},s_{in}-s_{2}^{\star})$.\\

The Jacobian matrix of dynamics (\ref{dyn2}) at the non-trivial
equilibrium
$(s_{1}^{\star},x_{1}^{\star},s_{2}^{\star},x_{2}^{\star})$
is of the following form in $(z_{1},z_{2},s_{1},s_{2})$ coordinates
\[
\left[\begin{array}{cc}
A_{s} & 0\\
\star & J^{\star}
\end{array}\right] \quad \mbox{with} \quad
J^{\star}=\left[\begin{array}{cc}s_{1}^{\star}-s_{in} & 0\\[2mm]
\frac{1}{1-r} & 2s_{2}^{\star}-\frac{1}{1-r}-s_{in}
\end{array}\right]
\]
Recall that $A_{s}$ is Hurwitz. The eigenvalues of $J^{\star}$ are
$-(s_{in}-s_{1}^{\star})$, $-(s_{in}-s_{2}^{\star})-(1/(1-r)-s_{2}^{\star})$ that are both negative numbers, accordingly to Proposition 1.
The exponential stability of the non-trivial equilibrium is thus proved. \qed

\subsection{The parallel configuration} 

\noindent {\bf Proposition  6.}
When $s_{in}>1$ and $d>0$, any trajectory of
(\ref{dyn2d}) with initial condition in $D$ such that 
$x_{1}(0)>0$ and $x_{2}(0)>0$
converges exponentially to the unique non-trivial steady-state 
$(s_{1}^{\star},x_{1}^{\star},s_{2}^{\star},x_{2}^{\star})$ given by Proposition 2.\\

{ \em Proof.}
Considering the time vector $z(\cdot)$,
the $(s_{1},s_{2})$ sub-system of dynamics (\ref{dyn2d}) can be written as solution of a non-autonomous planar dynamics
\begin{equation}
\label{nonautonomous2}
\left\{\begin{array}{lll}
\dot s_{1} & = & s_{1}(z_{1}(t)+s_{1}-s_{in})+\alpha_{1}(s_{in}-s_{1})+\frac{d}{r}(s_{2}-s_{1})\\
\dot s_{1} & = & s_{2}(z_{2}(t)+s_{2}-s_{in}) + \alpha_{2}(s_{in}-s_{2})+\frac{d}{1-r}(s_{1}-s_{2})
\end{array}\right.
\end{equation}
We know that $z$ converges to $0$ and consequently the vector
$S$ of variables $s_{1}$, $s_{2}$ converges to the set 
${\cal S}=[0,s_{in}]\times[0,s_{in}]$. We study now the limiting autonomous dynamics
\begin{equation}
\label{autonomous2}
\left\{\begin{array}{lll}
\dot s_{1} & = & (s_{in}-s_{1})(\alpha_{1}-s_{1})+\frac{d}{r}(s_{2}-s_{1})\\
\dot s_{2} & = & (s_{in}-s_{2})(\alpha_{2}-s_{2})+\frac{d}{1-r}(s_{1}-s_{2})
\end{array}\right.
\end{equation}
on the domain ${\cal S}$. 
Let ${\cal B}$ be the boundary 
$\{s_{1}=s_{in}\}\cup\{s_{2}=s_{in}\}$.
On the domain ${\cal S}\setminus{\cal B}$, we consider the vector $\sigma$ of variables $\sigma_{i}=\log(s_{in}-s_{i})$, whose dynamics can be written as follows
\begin{equation}
\label{dynsigma}
\dot\sigma = F(\sigma)=
\left[\begin{array}{c}
-\alpha_{1}+s_{in}-e^{\sigma_{1}}-\frac{d}{r}\left(1-e^{\sigma_{2}-\sigma_{1}}\right)\\
-\alpha_{2}+s_{in}-e^{\sigma_{2}}-\frac{d}{1-r}\left(1-e^{\sigma_{1}-\sigma_{2}}\right)
\end{array}\right]
\end{equation}
One can easily compute
\[
\mbox{div}(F)=-e^{\sigma_{1}}-e^{\sigma_{2}}-\frac{d}{r}e^{\sigma_{2}-\sigma_{1}}-\frac{d}{1-r}e^{\sigma_{1}-\sigma_{2}} <0
\]
From Poincar\'e-Bendixon theorem and Dulac criterion, we conclude that 
bounded trajectories of (\ref{dynsigma}) cannot have limit cycle or closed path and necessarily converge to an equilibrium point.
Consequently, any trajectory of (\ref{autonomous2}) in ${\cal S}$
either converges to the rest point 
$S^{\star}=(s_{1}^{\star},s_{2}^{\star})$ 
or approaches the boundary ${\cal B}$. Notice that one has
\[
s_{i}=s_{in}, \, s_{j}<s_{in} \; \Rightarrow \;
\dot s_{i}<0 \qquad (i\neq j)
\]
So the only possibility for approaching ${\cal B}$ is to converge to the other rest point $S^{0}=(s_{in},s_{in})$. 
This shows that the only 
non-empty, closed, connected, invariant and chain recurrent 
subsets of ${\cal S}$ are the singletons $\{S^{\star}\}$ and
 $\{S^{0}\}$.\\

\begin{center}
\begin{figure}[h]
\begin{center}
\includegraphics[height=5cm]{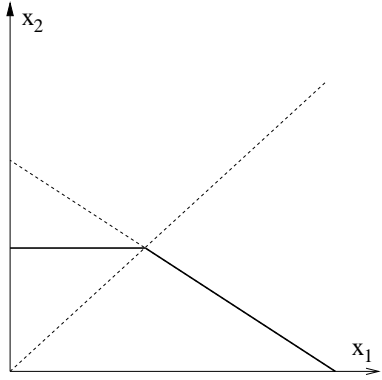}
\caption{Iso-value of the $V$ function.\label{fig-V}}
\end{center}
\end{figure}
\end{center}

Applying Theorem 4 we conclude that any trajectory of (\ref{nonautonomous2}), issued from initial condition of dynamics (\ref{dyn2d})
in ${\cal D}$, converges asymptotically to $S^{\star}$ or $S^{0}$.
Consider now any initial condition with $x_{1}(0)>0$ and  $x_{2}(0)>0$. 
We show that the solution $(s_{1}(\cdot),s_{2}(\cdot))$
of (\ref{dyn2d}) cannot converge to  $S^{0}$. 
If it is the case, there exists $T<+\infty$ such that one has
\[
s_{1}(t)>\alpha_{1} \quad \mbox{and} \quad rs_{1}(t)+(1-r)s_{2}(t)>1
\quad \mbox{for any } t\geq T
\]
under the assumption $s_{in}>1$.
Let us consider the function
\[
V(x_{1},x_{2})=\min(rx_{1}+(1-r)x_{2},x_{1})
\]
(see Figure \ref{fig-V})
and $v(t)=V(x_{1}(t),x_{2}(t))$ that is positive and tends to $0$ when $t$ tends to $+\infty$
If $x_{1}(t)< x_{2}(t)$, one has $v(t)=x_{1}(t)$ and
\[
\dot v=\dot x_{1}>(s_{1}(t)-\alpha_{1})x_{1}>0\mbox{ for } t\geq T
\]
If $x_{1}(t)>x_{2}(t)$, one has $v(t)=rx_{1}(t)+(1-r)x_{2}(t)$ and
\[
\begin{array}{lll}
\dot v & = & r\dot x_{1}+(1-r)\dot x_{2}
       = r(s_{1}-\alpha_{1})x_{1}+(1-r)(s_{2}-\alpha_{2})x_{2}\\
       & > & (rs_{1}+(1-r)s_{2}-1)x_{2} >0 \mbox{ for } t\geq T
\end{array}
\]
We conclude that the function $t\mapsto v(t)$ is non-decreasing for $t\geq T$ and consequently cannot converge to zero, thus a contradiction.\\

The Jacobian matrix of dynamics (\ref{dyn2d}) at the non-trivial
equilibrium
$(s_{1}^{\star},x_{1}^{\star},s_{2}^{\star},x_{2}^{\star})$
is of the following form in $(z_{1},z_{2},s_{1},s_{2})$ coordinates
\[
\left[\begin{array}{cc}
A_{p} & 0\\
\star & J^{\star}
\end{array}\right] \quad \mbox{with} \quad
J^{\star}=\left[\begin{array}{cc}-\frac{d}{r}\phi_{2}^{\prime}(s_{1}^{\star}) & \frac{d}{r}\\[2mm]
\frac{d}{1-r} & -\frac{d}{1-r}\phi_{1}^{\prime}(s_{2}^{\star})
\end{array}\right]
\]
Recall that $A_{p}$ is Hurwitz. One has
\[
\mbox{det}(J^{\star})=\frac{d^{2}}{r(1-r)}(\phi_{1}^{\prime}(s_{2}^{\star})\phi_{2}^{\prime}(s_{1}^{\star})-1) \quad \mbox{and} \quad
\mbox{tr}(J^{\star})=-\frac{d}{r}\phi_{2}^{\prime}(s_{1}^{\star})-\frac{d}{1-r}\phi_{1}^{\prime}(s_{2}^{\star}) \ .
\]
The function $\phi_{1}(\cdot)$ being concave, one has
$\phi_{1}(s_{in}) \leq \phi_{1}(s_{2}^{\star})+\phi_{1}^{\prime}(s_{2}^{\star})(s_{in}-s_{2}^{\star})$. Along with the inequalities $s_{in}>s_{2}^{\star}$ and $\phi_{1}(s_{in})=s_{in}>s_{1}^{\star}=\phi_{1}(s_{2}^{\star})$, one deduces $\phi_{1}^{\prime}(s_{2}^{\star})>0$. 
Recall from Corollary 2 that one has
$g^{\prime}(s_{1}^{\star})= \phi_{1}^{\prime}(s_{2}^{\star})\phi_{2}^{\prime}(s_{1}^{\star})-1> 0$.
Then the inequality $\phi_{2}^{\prime}(s_{1}^{\star})>0$ is
necessarily satisfied. Finally, we have shown $\mbox{det}(J^{\star})>0$ and $\mbox{tr}(J^{\star})<0$, that guarantee the
exponential stability of the non-trivial equilibrium 
$(s_{1}^{\star},x_{1}^{\star},s_{2}^{\star},x_{2}^{\star})$. \qed\\

\noindent {\em Remark.} The wash-out equilibrium $(s_{in},0,s_{in},0)$ is not necessarily hyperbolic. This explains why we cannot use the Convergence Theorem for asymptotically autonomous dynamics given in Appendix F of \cite{Smith}.

\end{document}